\documentclass[reqno]{amsart}
\usepackage{amsfonts,amsmath,amssymb,amsrefs,dsfont}
\usepackage{color}
\usepackage[utf8]{inputenc}  
\usepackage[T1]{fontenc}  
\usepackage{graphicx} 
\numberwithin{equation}{section}

\usepackage[colorinlistoftodos]{todonotes}

\newtheorem{thm}{Theorem}[section]
\newtheorem*{thm*}{Theorem}

\newtheorem*{cor*}{Corollary}

\begin{document}
\title[Sharp quantization for Lane-Emden problems in dimension two]
{Sharp quantization for Lane-Emden problems in dimension two}

\author{Pierre-Damien THIZY}
\address{Pierre-Damien Thizy\\
Univ Lyon, Universit\'e Claude Bernard Lyon 1, CNRS UMR 5208, Institut Camille Jordan, 43 blvd. du 11 novembre 1918, F-69622 Villeurbanne cedex, France}
\email{pierre-damien.thizy@univ-lyon1.fr}
\date{February 2018}

\begin{abstract}
 In this short note, we prove a sharp quantization for positive solutions of Lane-Emden problems in a bounded planar domain. This result has been conjectured by De Marchis, Ianni and Pacella \cite[Remark 1.2]{DeMarchisIanniPacella}.
\end{abstract}

\maketitle

\section*{Introduction}
Let $\Omega$ be an open, non-empty, {connected} and bounded subset of $\mathbb{R}^2$ with smooth boundary $\partial\Omega$ and let $\Delta=-(\partial_{xx}+\partial_{yy})$ be the (positive) laplacian. In this paper, we are interested in the asymptotic behavior as $p\to +\infty$ of a sequence $(u_p)_p$ of smooth functions, positive in $\Omega$, and satisfying the so-called Lane-Emden problem
\begin{equation}\label{LaneEmdenPbm}
\begin{cases}
&\Delta u_p = |u_p|^{p-1}u_p~~\text{ in }\Omega\,,\\
&u_p=0\quad\quad\quad\quad\quad\text{on }\partial\Omega\,,
\end{cases}
\end{equation} 
together with the bounded-energy type assumption 
\begin{equation}\label{BdEnergy}
p\int_\Omega |\nabla u_p|^2 ds=O(1)\,,
\end{equation}
for all $p$. Up to now, the most general results on this problem were obtained by De Marchis, Ianni and Pacella \cite{DeMarchisIanniPacella}. In particular, for such a given $(u_p)_p$ satisfying \eqref{LaneEmdenPbm}-\eqref{BdEnergy}, it is proved in \cite{DeMarchisIanniPacella} that, up to a subsequence, there exists an integer $n\ge 1$ and a subset $\mathcal{B}=\{x_1,...,x_n\}$ of $\Omega$ such that the following quantization
\begin{equation}\label{Quantization}
\lim_{p\to +\infty}p\int_\Omega |\nabla u_p|^2 ds=8\pi \sum_{j=1}^n m_j^2
\end{equation} 
holds true, where the $m_j$'s can be obtained through
\begin{equation}\label{DefMj}
m_j=\lim_{\beta\to 0^+} \lim_{p\to +\infty}\|u_p\|_{C^0(B_{x_j}(\beta))}\,,
\end{equation}
for $j\in\{1,...,n\}$, where $B_{x_j}(\beta)$ is the ball of center $x_j$ and radius $\beta$. Observe that, in particular, the $x_j$'s are not in $\partial\Omega$. It is also proved in \cite{DeMarchisIanniPacella} that we necessarily have that
\begin{equation}\label{MjGe}
m_j\ge \sqrt{e}\,,
\end{equation}
for all $j\in\{1,...,n\}$ and that
\begin{equation}\label{ConvLocWeak}
\lim_{p\to +\infty} u_p=0\text{ in }C^2_{loc}(\bar{\Omega}\backslash \mathcal{B})\,.
\end{equation}
 In \cite[Remark 1.2]{DeMarchisIanniPacella}, it is conjectured that we must have equality in \eqref{MjGe}, so that, in some sense, the constant $8\pi e$ plays here the same role as the Sobolev constant in dimensions greater than $2$ (see Struwe \cite{Struwe}). This is the point in the following theorem.
\begin{thm}\label{MainThm}
Let $\Omega$ be a smooth bounded domain of $\mathbb{R}^2$. Let $(u_p)_p$ be a sequence of smooth functions positive in $\Omega$, and satisfying \eqref{LaneEmdenPbm} and \eqref{BdEnergy}. Then, up to a subsequence, there exists an integer $n\ge 1$ such that
\begin{equation}\label{QuantizationAmel}
\lim_{p\to +\infty}p\int_\Omega |\nabla u_p|^2 ds=(8\pi e)\times n\,.
\end{equation}
Moreover, there exists a subset $\{x_1,...,x_n\}$ of $\Omega$ such that 
\begin{equation}\label{MjEq}
m_j=\sqrt{e}\,,
\end{equation}
where $m_j$ is given by \eqref{DefMj}, for all $j\in\{1,...,n\}$.
\end{thm}
In addition, by \cite{DeMarchisIanniPacella}, we get from \eqref{MjEq} that
\begin{equation}\label{ConvLoc}
\lim_{p\to +\infty} p u_p= 8\pi \sqrt{e} \sum_{j=1}^n \mathcal{G}_{x_j} \text{ in }C^2_{loc}(\bar{\Omega}\backslash \mathcal{B})\,,
\end{equation}
and that
\begin{equation}\label{Poho}
\nabla_{x_j}\left(\mathcal{H}_{x_j}(x_j)+\sum_{i\neq j}\mathcal{G}_{x_i}(x_j) \right)=0\,,
\end{equation}
for all $j\in\{1,...,n\}$, where $\mathcal{G}$ is the Green's function of $\Delta$ with Dirichlet boundary conditions and where $\mathcal{H}$ is its regular part, which is smooth in $\Omega^2$ and given by
$$\mathcal{G}_x(y)=\frac{1}{2\pi} \log\frac{1}{|x-y|}+\mathcal{H}_x(y)\,, $$
for all $x\neq y$.\\
\indent Concerning the previous works, Ren and Wei \cite{RenWei} and \cite{RenWei2} where able to prove that \eqref{QuantizationAmel} with $n=1$ holds true if the $u_p$'s are minimizers, i.e. if we assume in addition that $u_p$ is proportional to a solution of the problem $$\min_{\{v\in H^1_0\text{ s.t. }\int_\Omega v^p ds =1\}} \int_\Omega |\nabla v|^2 ds .$$
Answering to a former question, Adimurthi and Grossi \cite{AdimGrossi} were able to prove that
\begin{equation}\label{LInftyNorm}
\lim_{p\to +\infty}\|u_p\|_{C^0(\Omega)}=\sqrt{e}\,,
\end{equation}
in the case of minimizers, while they discovered the way to perform the first rescaling for the $u_p$'s as $p\to +\infty$, and the key link with the Liouville equation. Observe that \eqref{DefMj}, \eqref{ConvLocWeak} and \eqref{MjEq} clearly imply \eqref{LInftyNorm} in general case.  Now in the radial case where $\Omega$ is a disk, observe that the $u_p$'s are necessarily minimizers, since \eqref{LaneEmdenPbm} admits only one solution (see Gidas-Ni-Niremberg \cite{GidasNiNiremb} and the nice survey by Pacella \cite{PacellaUniqueness}); according to the previous discussion, we necessarily then have that $n=1$ in \eqref{QuantizationAmel}. In contrast, if $\Omega$ is not simply connected, Esposito, Musso and Pistoia \cite{EspMussPisConcentSolLargeExpo} were able to prove that, for all given integer $n\ge 1$, there exists a sequence of positive functions $(u_p)_p$ satisfying \eqref{LaneEmdenPbm}-\eqref{BdEnergy} such that \eqref{QuantizationAmel} holds true, together with \eqref{MjEq}-\eqref{LInftyNorm}. Thus, in some sense, Theorem \ref{MainThm} is sharp. We mention that very interesting complementary results were obtained recently by Kamburov and Sirakov \cite{KamSir}. At last, we also mention that, even if the situation is far from being as well understood in the nodal case, where we no longer assume that the $u_p$'s are positive, some asymptotic-analysis \cite{AsymptoticAnaSignChanDMIP,LEPbmAsymptBehLowEnerGrossiGrumiauPacella}, as well as some constructive \cite{MorseIndexDMIP,EspMussPisConcentSolLargeExpoNodal,GladialiIanni} results were obtained. 
\\
\indent To conclude, as explained in De Marchis, Ianni and Pacella \cite{DeMarchisIanniPacellaSurvey}, the techniques to get the quantization result in \cite{DeMarchisIanniPacella} are not without similarity with the ones developed by Druet \cite{DruetDuke} to get the analogue quantization for 2D Moser-Trudinger critical problems. Both results \cite{DeMarchisIanniPacella,DruetDuke} can be improved by showing that all the blow-up points necessarily carry the minimal energy. It is done here in the Lane-Emden case and in Druet and Thizy \cite{DruThiI} in the Moser-Trudinger case. Unfortunately, the authors of \cite{DruThiI} were not able to find an as easy argument as here, in the more tricky Moser-Trudinger critical case.

\section{Proof of Theorem \ref{MainThm}}
Let $(u_p)_p$ be a sequence of smooth functions, positive in $\Omega$ and satisfying \eqref{LaneEmdenPbm}-\eqref{BdEnergy}. Then by \cite{DeMarchisIanniPacella}, \eqref{Quantization}-\eqref{ConvLocWeak} hold true. Thus, the proof of \eqref{QuantizationAmel}-\eqref{MjEq}, i.e. that of Theorem \ref{MainThm}, reduces to the proof of
\begin{equation}\label{LInftyLe}
\lim_{p\to +\infty}\|u_p\|_{C^0(\Omega)}\le \sqrt{e}\,.
\end{equation} 
Here and in the sequel, we argue up to a subsequence. Now, let $(y_p)_p$ be a sequence in $\Omega$ such that $u_p(y_p)=\|u_p\|_{C^0(\Omega)}$, for all $p$. By \eqref{DefMj}-\eqref{ConvLocWeak}, we have that
\begin{equation}\label{NotToBdry}
\lim_{p\to +\infty} d(y_p,\partial \Omega):=2\delta_0>0\,,
\end{equation}
where $d(y,\partial \Omega)$ denotes the distance from $y$ to $\partial\Omega$. Now, let $\mu_p>0$ be given by 
\begin{equation}\label{Mup}
\mu_p^2 ~p~ u_p(y_p)^{p-1}=8\,.
\end{equation}
By \eqref{DefMj} and \eqref{MjGe}, we get from \eqref{Mup} that
\begin{equation}\label{Mup2}
\log\frac{1}{\mu_p^2}=p \log u_p(y_p)\left(1+O\left(\frac{\log p}{p} \right)\right)\,,
\end{equation}
and in particular, that $\mu_p\to 0$ as $p\to +\infty$. Let $\tau_p$ be given by
$$u_p(y_p+\mu_p y)=u_p(y_p)\left(1-\frac{2\tau_p(y)}{p} \right)\,, $$
so that 
\begin{equation}\label{CondConvTau}
\tau_p\ge 0\text{ and }\tau_p(0)=0\,,
\end{equation}
 by definition of $(y_p)_p$. By \eqref{LaneEmdenPbm} and \eqref{Mup}, we have that
\begin{equation}\label{EqTau}
\Delta(-\tau_p)=4\left(1-\frac{2\tau_p}{p} \right)^p\text{ in }\Omega_p:=\frac{\Omega-y_p}{\mu_p}\,,
\end{equation}
so that, by \eqref{CondConvTau}, positivity of $u_p$ and concavity of the $\log$ function, we get that
\begin{equation}\label{EstimLaplTau}
0<\Delta(-\tau_p)\le 4\,.
\end{equation}
By \eqref{NotToBdry}, \eqref{CondConvTau}-\eqref{EstimLaplTau} and standard elliptic theory, including the Harnack principle, we get that there exists a function $\tau_\infty\in C^2(\mathbb{R}^2)$ such that
\begin{equation}\label{ConvTau}
\lim_{p\to+\infty} \tau_p=\tau_\infty\text{ in }C^2_{loc}(\mathbb{R}^2)\,,
\end{equation}
and then that
\begin{equation}\label{CondEqTau}
\Delta(-\tau_\infty)=4\exp(-2\tau_\infty)\text{ in }\mathbb{R}^2\,,\quad \tau_\infty(0)=0\,,\quad \nabla \tau_\infty(0)=0\,,
\end{equation}
using also that $\nabla \tau_p(0)=0$, by definition of $(y_p)_p$. Let $R>0$ be given. Integrating by parts, using \eqref{Mup}, \eqref{ConvTau}, \eqref{BdEnergy}, $\Delta u_p\ge 0$ and $u_p\ge 0$, we get that
\begin{equation}
\begin{split}
\liminf_{p\to +\infty}(8u_p(y_p)^2)\int_{B_0(R)} 4\exp(-2\tau_\infty) dy &\le\lim_{p\to +\infty}p\int_{B_{y_p}(R\mu_p)}(\Delta u_p) u_p dy\,,\\
&\le \lim_{p\to +\infty}p\int_{\Omega} |\nabla u_p|^2 dy<+\infty\,. 
\end{split} 
\end{equation}
 By using that $u_p(y_p)^2\ge (1+o(1))e$~, and by observing that the above RHS does not depend on $R>0$, which can be arbitrarily large, we get that
\begin{equation}\label{ConfVolFinite}
\int_{\mathbb{R}^2} \exp(-2\tau_\infty) dy<+\infty \,.
\end{equation}
By Chen and Li \cite{ChenLi}, \eqref{CondEqTau} and \eqref{ConfVolFinite} imply that 
\begin{equation}\label{TauInftyFormula}
\tau_\infty=\log(1+|\cdot|^2)\,.
\end{equation}
 Then, we let $t_p$ be given by $$t_p(y)=\log\left(1+\frac{|y-y_p|^2}{\mu_p^2} \right)\,.$$
From now on, if $f$ is a given continuous function in $\Omega$, we let $\bar{f}$ be the unique continuous function in $[0,d(y_p,\partial\Omega))$ given by
$$\bar{f}(r)=\frac{1}{2\pi r} \int_{\partial B_{y_p}(r)} f d\sigma\,, \text{ for all }r\in(0,,d(y_p,\partial\Omega)) \,.$$
Let $(r_p)_p$ be any sequence such that $r_p\in [0,d(y_p,\partial\Omega))$ for all $p$. By \eqref{LaneEmdenPbm}, \eqref{Mup}, \eqref{ConvTau}, \eqref{TauInftyFormula} and Fatou's lemma, we get that
\begin{equation}\label{ComputIntegLapl}
\begin{split}
&-2\pi r_p \frac{d \bar{u}_p}{dr}(r_p) \\
&\quad\quad\quad\quad~=\int_{B_{y_p}(r_p)} \overline{(\Delta u_p)}~ 2\pi r dr\,,\\
&\quad\quad\quad\quad~\ge \frac{2 u_p(y_p)}{p}\left(\int_0^{r_p/\mu_p} \frac{8\pi r~ dr}{(1+r^2)^2}+o\left(\frac{(r_p/\mu_p)^2}{1+(r_p/\mu_p)^2} \right)\right)\,,
\end{split}
\end{equation}
using that the laplacian commutes with the average in spheres. Then, using the fundamental theorem of calculus and $\bar{u}_p(0)=u_p(y_p)$, we easily get from \eqref{ComputIntegLapl} that
\begin{equation}\label{intermsupp}
\bar{u}_p(r_p)\le u_p(y_p)\left(1-\frac{2 \bar{t}_p(r_p)}{p}+o\left(\frac{\bar{t}_p(r_p)}{p} \right)\right)\,.
\end{equation} 
Picking now $r_p=\delta_0$ for all $p$, according to \eqref{NotToBdry}, we get from \eqref{intermsupp} that
\begin{equation}\label{CCl}
\bar{u}_p(\delta_0)\le {u_p(y_p)}\left(1-2\log u_p(y_p)\left(1+o(1)\right)+O\left(\frac{1}{p}\right) \right)\,, 
\end{equation}
by \eqref{Mup2}, writing merely $\log\frac{1}{\mu_p^2}+\log\left(\mu_p^2+\delta_0^2 \right)=\bar{t}_p(\delta_0)$. Since $\bar{u}_p\ge 0$ and since $u_p(y_p)=\|u_p\|_{C^0(\Omega)}$, we easily get from \eqref{CCl} that \eqref{LInftyLe} holds true. As explained at the beginning of the proof, this concludes the proof of Theorem \ref{MainThm}.


\end{document}